\documentclass[12pt,oneside]{article}
\newcommand{\Path}{}
\usepackage{ \Path Paper_style_2}
\usepackage{ \Path Keywords_WB}
\usepackage{ \Path Environments}
\newcommand{\figs}{}
\newcommand{\diags}{}
\usepackage{tikz-cd}
\usepackage{ \Path tikzcd2}
\usepackage{sidecap}

\DeclareMathOperator{\Cylndr}{Cyl}
\DeclareMathOperator{\Pipe}{Pipe}

\newcommand{\law}{\mathfrak{m}}
\newcommand{\transP}{\mathfrak{p}}
\newcommand{\Paths}{ \text{Paths} }
\newcommand{\prob}[2][\text{prob}]{ #1 \left( #2 \right)}
\newcommand{\Condprob}[3][\text{prob}]{ #1 \left( #2 \,\rvert\, #3 \right)}
\newcommand{\generatorF}{\Gamma_{\text{shift}, 1}}
\newcommand{\generator}{ \left( \Omega_{\text{shift}}, \generatorF^t, \nu_{\text{shift}} \right) }
\begin{document}
\title{Discrete-time dynamics, step-skew products, and pipe-flows}
\author{Suddhasattwa Das\footnotemark[1]}
\footnotetext[1]{Department of Mathematics and Statistics, Texas Tech University, Texas, USA}
\date{\today}
\maketitle
\begin{abstract} 
	Dynamical processes can be classified in various ways as deterministic or stochastic, and continuous or discrete time. All these types can be studied by the path-spaces they generate, and stationary measures on that path-space. Such measures are called the law of the dynamics. This article presents how a general ergodic dynamical system may be approximated in terms of their law, by a simple and restricted family of deterministic continuous-time skew-product systems. In these systems, a deterministic, mixing flow intermittently drives a deterministic flow through a topological space created by gluing cylinders. The resulting orbits mimic the law of the original dynamics. This comparison is made possible by introducing a secondary intermediary approximation of the ergodic dynamics. This third system is a step-skew dynamical system, in which a finite state Markov process drives a dynamics on topological disk. Each of these three representations have their advantages. It is proved that the distribution induced on the space of paths by these three dynamics can be made arbitrarily close to each other. This analysis reconfirms the old principle that it is impossible to decide whether a general timeseries is generated by a deterministic or stochastic process, and is of continuous or discrete time.
\end{abstract}


\section{Introduction} \label{sec:intro}

\begin{figure}[!t]
	\center
	\begin{tikzpicture}[scale=0.55, transform shape]
		\input{\figs Pipe_flow_outline.tex}
	\end{tikzpicture}
	\caption{Alternative descriptions of ergodic dynamics. The starting point of this analysis is an ergodic system, as defined in Assumption \ref{A:f}. This is a general means of describing most deterministic physical phenomenon. One of its key components is the invariant measure $\mu$. It has two significances -- its support represents the phase space of the dynamics under observation. Its distribution determines the statistical properties of its generated data, as well as its dynamic complexity. The article presents how the dynamics can be approximated by the other two types of dynamics, labelled II and III. Type II is a skew product system in which an autonomous finite-state Markov process drives a dynamical system on $d$-dimensional disks. Type II performs a dual stochastic and topological approximation of Type I. Type III is a deterministic flow through topological spaces called cells and cylinders. A trajectory is observed when it exits these cells. The timeseries created by the series of exit points provide a statistical approximation of the Type II. The meaning of these approximations is explained in the smaller white boxes. They connect various secondary characteristic of these dynamical systems.}
	\label{fig:outline}
\end{figure} 

Theoretical studies of Dynamical systems have lead to the invention of various simple templates or toy-models, which reveal some of the many complex behavior of dynamical systems.  On the other hand numerical methods, especially data-driven ones, focus on the statistical and topological properties emerging from timeseries generated by measurements of a dynamical system. Their focus is convergence and consistency, and overcoming challenges such as under-sampling and high dimensionality. There is a major gap in these theoretical and numerical approaches. The theoretical templates are too specific and formula-dependent to be found in natural or physical systems. On the other hand the reconstruction provided by numerical methods mostly fail to include the rich dynamical properties discovered in theory, and only focus on recreating the dynamical law. The present work aims to bridge this gap by showing that any ergodic, discrete time dynamical system may be mimicked by samples of a continuous-time deterministic dynamical, in a manner which approximates both its \textit{law} as well as the topology of the invariant set. The approximant dynamics will be shown to be highly interpretable and conducive to data-driven methods. We shall assume throughout that :

\begin{Assumption} \label{A:f}
	There is a continuous map $f: \real^d \to \real^d$, an invariant ergodic measure $\mu$ of $f$ whose support $X$ is compact.
\end{Assumption}

Two of the most important features of an ergodic system as above are the topology of the set $X$, and the \textit{law} it generates in the ambient space $\real^d$. The law of any dynamic process, explained in detail in Section \ref{sec:main}, describes the statistical properties of the timeseries it generates. As described in \cite{BerryDas2023learning, BerryDas2024review} and the references therein, one of the biggest challenges in the approximation theory for dynamics is in achieving a simultaneous approximation of the dynamics transition rule, the law, and the topology. The ideas presented in this article provides an explicit constructive means of filling this gap. Figure \ref{fig:outline} presents an outline of the proposed approximation scheme.

The approximation proceeds in two steps. In an intermediate step, the dynamics in Assumption \ref{A:f} is approximated as a \emph{step-skew product} 
\begin{equation} \label{eqn:Mrkv:1}
	\begin{aligned}
		s_{n+1} &\sim \transP(s_n) \\
		y_{n+1} &= \phi_{ s_n\to s_{n+1} } ( y_n )
	\end{aligned} ,\quad \quad 
	(s_n, y_n) \in \calS \times \calD.
\end{equation}
This is a skew product system in which the first set of coordinates (namely $s$) evolves autonomously, and drives the second set of coordinates. The variable $s$ is drawn from a finite state-space, and undergoes a Markov walk via a Markov transition matrix $\transP$. The second variable $y$ is drawn from some contractible space $\calD$ and undergoes deterministic transitions depending upon the Markov transition taking place in the $s$-variable. It was shown in \cite{Das2024zero} that the system \eqref{eqn:Mrkv:1} can be designed so that it simultaneously approximates the dynamics rule $f$, the targeted invariant region $X$, as well as the law induced by $f$. 

Thus the approximation strategy we propose first converts a targeted deterministic process on $\real^d$ into a process that is discrete-time stochastic and has two parts. The first part is an autonomous finite-state discrete-time Markov process. This process drives the second part, which is a dynamical system on $\calD$. The space $\calD$ is usually chosen to be $\real^d$ or the disk $\calD^d$. This leads  to a \textit{skew-product system} on $\calD^d$. We recall the following key result :

\begin{corollary} \label{corr:1}
	\cite{Das2024zero} Any ergodic dynamical system as in Assumption \ref{A:f} is the zero-noise limit of a Markov process on $\calD^d$, presented in format \eqref{eqn:Mrkv:1}. 
\end{corollary}

The approach that we take is distinguished  from other analyses by the fact that the limiting dynamics is not derived from the sequence of Markov processes. Rather, we take an arbitrary ergodic dynamics as in Assumption \ref{A:f}, and prescribe a sequence of Markov processes whose zero noise limit is $f:X\to X$. Each Markov process is a step-skew dynamical system, in which a finite-state Markov process drives a dynamical system on the $d$-dimensional disk $\calD^d$. The second stage and the main idea of this article is to obtain a continuous-time, deterministic realization of any such Markov process. Let $\Endo(\real^D)$ denote the group of $C^1$ endomorphisms of the Euclidean space $\real^D$. Suppose there is a driving dynamical system $\Gamma^t : \Omega \to \Omega$. Also suppose that there is a $\Endo(\real^D)$-valued cocycle $G$ over $\Gamma^t$. This means a map
\[ G : \real_0^+ \times \Omega \to \Endo(\real^D) , \]
such that
\[ G \paran{ t+s, \omega_0 } = G \paran{ t, 0, \Gamma^{s} \omega_0 } G \paran{ s, \omega_0 } , \quad \forall t,s \geq 0 , \; \forall \omega_0 \in \Omega . \]
This leads to a skew-product dynamics 
\begin{equation} \label{eqn:pipe_flow}
	\begin{tikzcd}
		\Omega \times \real^{D} \arrow[d, "\Psi_{T}^t"'] \\ 
		\Omega \times \real^{D}
	\end{tikzcd} ,\quad \quad 
	\left\{
	\begin{split}
		\omega(t) &= \Gamma^{Tt}(\omega_0) \\
		y(t) &= G(t, \omega_0) y_0
	\end{split}
	\right.
\end{equation}
The $T$ above represents a time-scaling parameter constant. Its purpose is to speed up the flow of the driving dynamics in $\omega$. We shall prove that

\begin{corollary} \label{corr:2}
	Any step-skew dynamical system as in \eqref{eqn:Mrkv:1} driving a dynamics on a contractible space $\calD$ can be weakly, conditionally approximated in law by a class of continuous-time, deterministic flows $\Psi_{T}^t$ of the form \eqref{eqn:pipe_flow}.
\end{corollary}

The concept of \emph{weak conditional convergence in law} is a mode of statistical approximation explained in detail in Section \ref{sec:random}. The approximation scheme that we present aims to approximate a discrete-time stochastic system via a continuous time, deterministic process. This scheme relies on the property of \emph{mixing}. If a deterministic system is mixing, then due to a phenomenon called decay of correlations, measurements of the same signal appear uncorrelated and random after a passage of time. Thus a mixing system can mimic the role of a stochastic system, while maintaining determinism and continuity. If one observes a continuous time flow at regular time-intervals of $\Delta t$, one obtains a discrete-time system, and its generator is called the \emph{flow-map} at time $\Delta t$. Our results can be summarized as follows : 

\begin{corollary} \label{corr:3}
	Consider the dynamical system $(X, f, \mu)$ in Assumption \ref{A:f}. Fix a spatial error-limit $\delta>0$, an uncertainty limit $\epsilon>0$, and a time $N\in \num$. Then there is a mixing dynamical system $\paran{\Omega, \Gamma^t, \nu}$ , a pipe flow as in \eqref{eqn:pipe_flow} on a topological space $\tilde{X}$ driven by $\Gamma^t$, an injective map $h:X\to \tilde{X}$ and a $T_0>0$ such that
	\begin{equation} \label{eqn:corr:3}
		\mbox{for } \mu- \mbox{a.e. } x_0 \,:\; \forall T>T_0, \quad \nu\SetDef{ \omega_0 }{ \abs{ h\paran{ f^n(x_0) } - y(3n) } < \delta, \; 1\leq n \leq N . } > 1-\epsilon .
	\end{equation}
	Here $\omega(t), y(t)$ is the trajectory of the flow \eqref{eqn:pipe_flow} starting from $\paran{\omega_0, h(x_0)}$. In other words, with probability at least $1-\epsilon$ the time-3 points on the orbit of the pipe-flow $\delta$-shadows the orbits of $f$ for $N$ time steps.
\end{corollary}

\begin{figure}[!t]
	\center
	\begin{tikzpicture}[scale=0.7, transform shape]
		\input{\diags Four_processes.tex}
	\end{tikzpicture}
	\caption{The four dynamical processes and their mutual approximations. Dynamical systems, i.e. process in which the future states are an outcome of the present state, may be broadly divided into four types, as shown in green text. There has been many interesting connections established between these processes. The connections are in the form of a dynamics of a certain type approximating another type, in terms of the invariant measure induced on the space of paths. This measure is known as the \emph{law} of the process. The diagram places the present work in relation to various similar discoveries. It relies on a sub-class of dynamics presented in \cite{Das2024zero} called \emph{step-skew} products. }
	\label{fig:four_process}
\end{figure} 

Corollary \ref{corr:3} is an easy consequence of our main result Theorem \ref{thm:3} presented later. Corollary \ref{corr:3} and Theorem \ref{thm:3} requires the construction of a special continuous-time flow whose descriptions spans Sections \ref{sec:random}, \ref{sec:junction} and \ref{sec:network}. Corollaries \ref{corr:1} and \ref{corr:3} indicate that it is impossible to decide whether a general timeseries is generated by a deterministic or stochastic process, and is of continuous or discrete time.

A proper approximation of a dynamical system requires a deeper understanding of the complexities that could occur and coexist, especially in chaotic dynamical systems \cite[e.g.]{KatokHassel1997, newhouse2004new, ACW2016, DasJim2017chaos, DasYorke2020}. This class of dynamical systems are characterized by infinitely many periodic orbits which co-exist in a dense formation. The phase space can be split everywhere into stable and unstable directions. Typical orbits get pulled towards either of these cycles for varying periods of time, before getting scattered by the unstable direction. A topological consequence is that orbits are seen to have almost periodic behavior for extended periods of time. An important statistical consequence is that statistics along such a typical trajectory appear uncorrelated, in spite of being obtained from a deterministic orbit. A naive approximation of the dynamics law ignores these rich details and just recreates the dynamics law.

This project is part of a broader effort \cite[e.g.]{froyland2001extract, Mischaikow2002topo, BoczkoKaliesMischaikow2007, FroylandPadberg09, kaczynski2016towards, BatkoEtAl2020linking, MrozekWanner2021creating, DGGS2024sec} to find alternative ways to describe a dynamical system, instead of the dynamics law alone. Dynamic processes can be broadly divided in many ways, such as discrete or continuous time, and deterministic or stochastic. This broad classification leads to four different types of processes. There has been several discoveries on how the timeseries generated by these four types can be approximated by one another, as summarized in Figure \ref{fig:four_process}. 

\begin{figure}[!t]
	\center
	\includegraphics[width=0.55\linewidth, height=0.6\textheight, keepaspectratio]{\figs 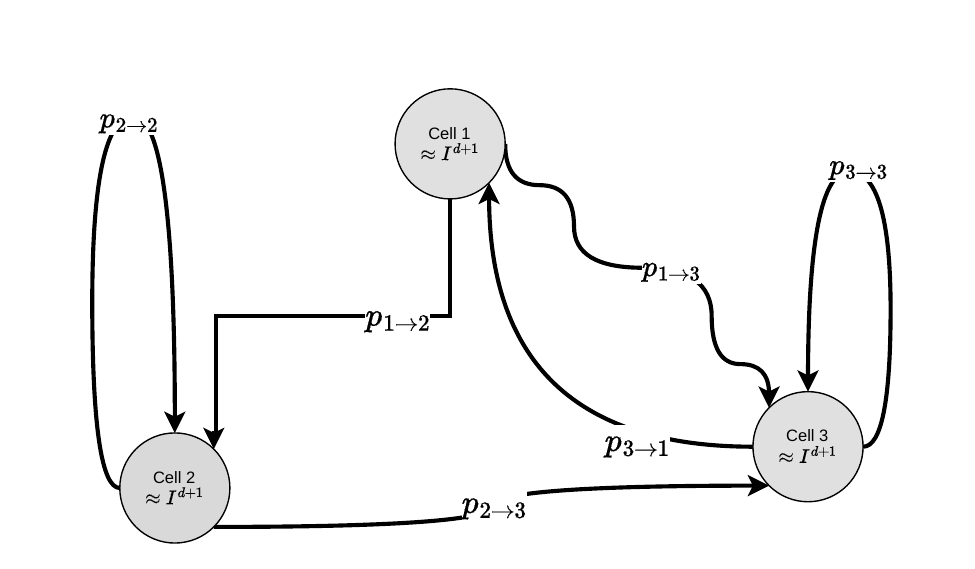}
	\includegraphics[width=0.35\linewidth, height=0.6\textheight, keepaspectratio]{\figs 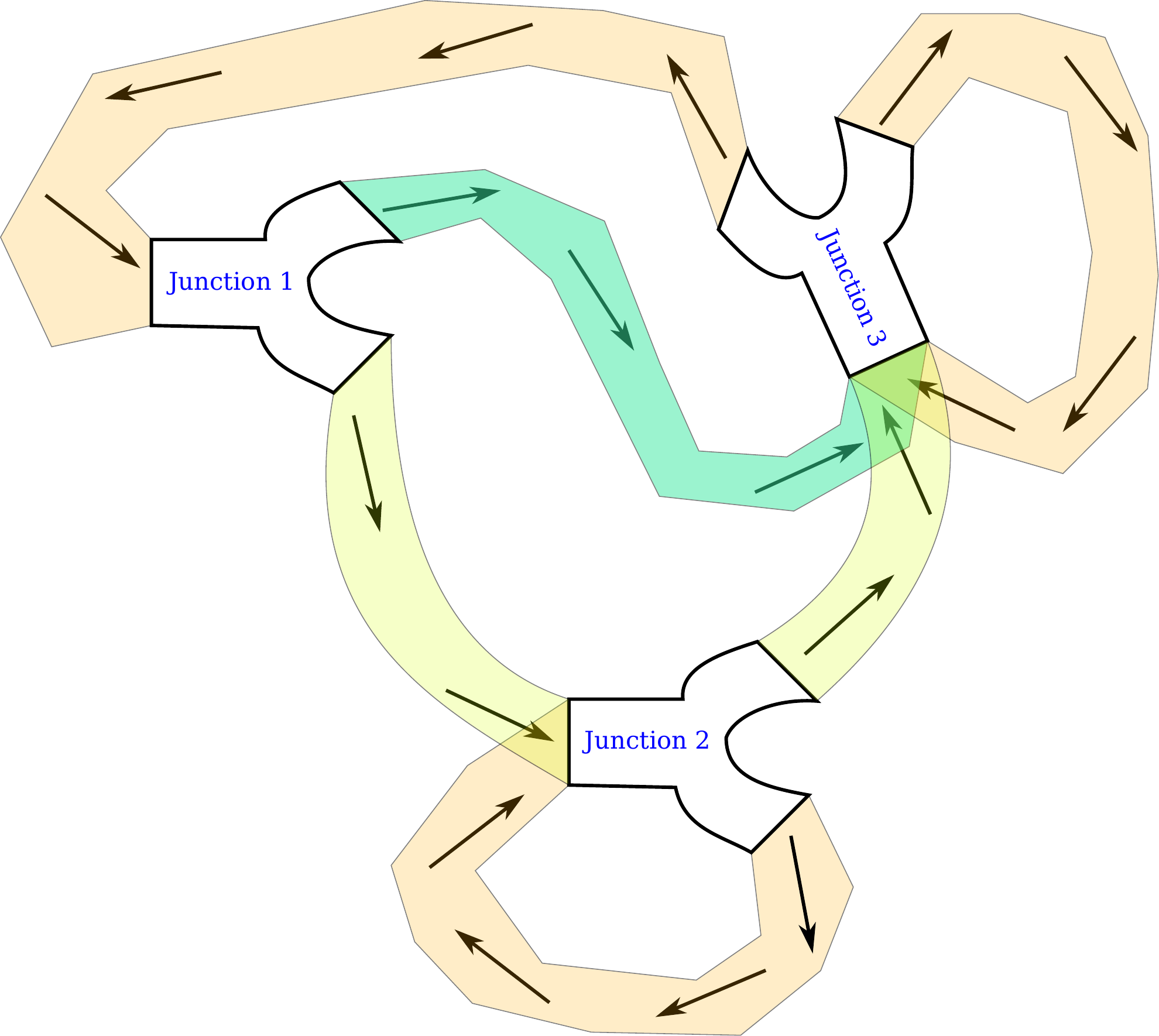}
	\caption{A transition network and its pipe-flow realization. A Markov process with three states have been depicted. The outgoing arrows from the $j$-th vertex have labels $\nu^{(j)}_i$ for $i\in 1,2,\ldots$. Each quantity $\nu^{(j)}_i$ represents the probability of transition to state-$i$ from state-$j$. The collection of transition probabilities create a vector $\nu^{(j)}$ whose entries sum to $1$. The idea of the construction is convert each of the states into a cell, which is topologically a $d+1$-dimensional cube. And each connection between cells is to made into a mapping torus corresponding to the maps $\phi_{j\to i}$. The end-faces of these mapping torus embed into $d$-dimensional cubes. See Section \ref{sec:network} for the details of the construction and the flow on the branched manifold shown on the right.}
	\label{fig:network}
\end{figure}

The ideas presented in the construction are the amalgamation of three broad ideas in dynamical systems theory. The first is the result from \cite{Das2024zero} that arbitrary ergodic systems are the zero noise limits of Markov processes. The second is the idea that chaotic systems can mimic various stochastic phenomenon. This has been realized in various formats and contexts, see \cite{Rudolph1976two, AlpernPrasad1989coding, AlpernPrasad2005towers, kieffer1980coding} for some notable examples.

The third key idea is directly involved in converting the step-skew product \eqref{eqn:Mrkv:1} into a deterministic flow. A finite state Markov walk involves random switches at every step. These switches, which are supposed to be independent random events, are realized as topological objects resembling a junction of pipes. The random switch is actuated by an external mixing flow, as indicated in \eqref{eqn:pipe_flow}. The higher the time-scaling factor $T$ is, the closer is the approximation to true randomness. This idea of using an external, fast, mixing system as a randomness generator finds similarity in other results \cite[e.g.]{MelbourneNicol2005almost, BurovEtAl2021kernel, MelbourneStuart2011skew}.

Our construction provides a tangible means of recreating the topology of an ergodic set $X$, along with orbits under the dynamics $f$.  We present several practical applications of this construction in Section \ref{sec:conclus}.

\paragraph{Outline} We begin by reviewing how dynamical systems have a natural interpretation as a step-skew dynamical system. We discuss their approximation properties in Section \ref{sec:step}. Having established the importance of step-skew products, we proceed to describe a continuous realization of step-skew products in Section \ref{sec:main}. This is where our main technical result Theorem  \ref{thm:3} is stated and proved. A main theoretical tool of this work is the behavior of mixing systems as pseudo-random generators. This is reviewed in Section \ref{sec:random}. Sections \ref{sec:junction} and Section \ref{sec:network} present the details of the construction. We end with some discussion on the potential applications of Theorem \ref{thm:3} in Section \ref{sec:conclus}. Some technical proofs are postponed to Section \ref{sec:appendix}.

\paragraph{Notations} Throughout the paper we shall use $\calD^d$ to denote the $d$-dimensional open disk. We use $I$ to denote the closed one dimensional interval $[0,1]$.

\section{Step-product dynamical systems} \label{sec:step}

\begin{figure}[!t]
	\center
	\includegraphics[width=0.95\linewidth, height=0.5\textheight, keepaspectratio]{\figs 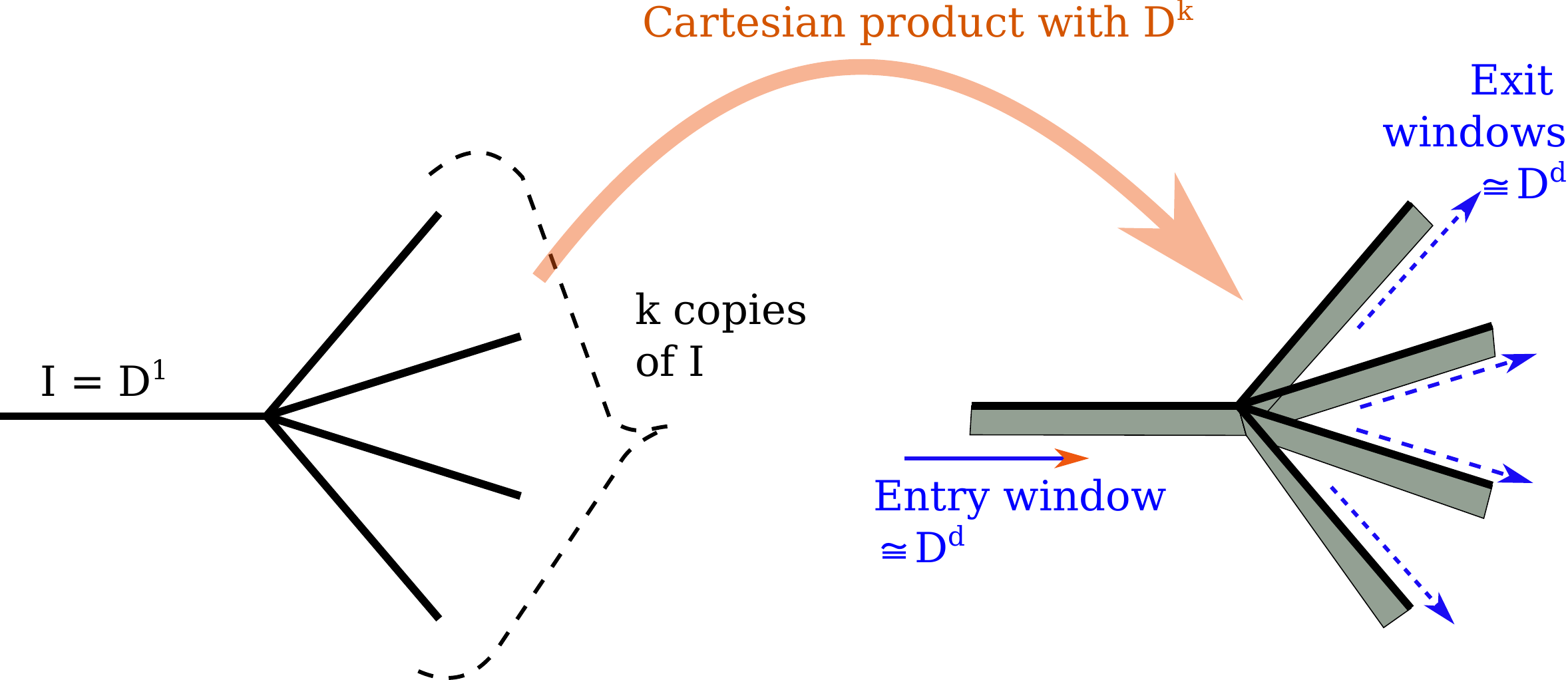}
	\caption{Construction of a junction. A $k$-junction in $d$-dimensions is a gluing of $k+1$ $d$-dimensional disks, as shown in the figure. See Section \ref{sec:junction} for a description of its use. Any such junction represents a state $s$ of the skew-product system \eqref{eqn:Mrkv:1}, which has $k$ possible outgoing states. The left most window is interpreted as the entry point, and the other $k$ terminals are interpreted as exit windows. See Figure \ref{fig:frontal}~(a) for the construction for general $k$. This entire topological space can be embedded in $\real^{d+3}$, as explained in the text.}
	\label{fig:cell}
\end{figure}

Throughout this section we assume Assumption \ref{A:f}.  We shall also assume 

\begin{Assumption} \label{A:cover}
	There is a finite measurable partition $\calU = \SetDef{U_i}{ 1\leq i \leq m}$ for $X$.
\end{Assumption}

Partitions provide a coarse graining of the phase space or invariant space. The space $X$ which is usually a continuum, is approximated by a finite set of cells. By keeping track of the transitions between the cells one obtains an outer approximation of the dynamics. For each cell $j$, we can compute transition probability vectors $\transP(j)$ in the following manner :
\begin{equation} \label{eqn:def:nu}
	\transP(j) \in \real^m, \quad \paran{ \transP(j) }_i := \mu \paran{ f(x) \in U_i \,\rvert\, x\in U_j } = \mu \paran{ f^{-1}(U_i) } / \mu \paran{ U_j } .
\end{equation}
Note that by design each vector $\transP(j)$ is a probability vector on the $\calS := \braces{1, \ldots, m}$. Our construction shall contain as a subsystem a discrete state Markov process on $\calS$. The power-set of $\calS$ is its assigned sigma-algebra. Consider the $m\times m$ matrix $\mathbb{P}$ whose $j$-th column is $\transP(j)$. The matrix $\mathbb{P}$ converts any probability measure $\beta$ on $\calS$ into another probability measure. Thus $\mathbb{P}$ plays the role of a Markov transition function on $\calS$. For each $j\in \calS$ we shall denote by $\transP(j)$ the (discrete) probability measure $\transP(j)$.

Consider the following sets for each transition $j\to i$ under $\mathbb{P}$ :
\begin{equation} \label{eqn:def:alpha_ij}
	\calX_{j\to i} := U_j \cap f^{-1} \paran{U_i} .
\end{equation}
These sets $\calX_{j\to i}$ will be used to  define maps $\phi_{j\to i}$ as follows 
\begin{equation} \label{eqn:def:ph_ji}
	\phi_{j\to i} : \real^d \to U_i, \quad \phi_{j\to i} \paran{ x } = f(x), \quad \mu-\mbox{a.e. } \, x\in \calX_{j\to i} .
\end{equation}
The functions $\phi_{j\to i}$ has two important features. Firstly, its range is confined to $U_i$. Secondly it agrees with the original map on a subset of the domain. This leads to a Markov process \eqref{eqn:Mrkv:1} on the product space $\calS \times \real^d$.  

\paragraph{Zero-noise limits} The notion of stability is made more precise using Arnold's paradigm \cite{Arnold1991random}. Let $\paran{ \Omega, \Sigma }$ be a measurable space, and $\tau : \Omega \to \Omega$ be a measurable map. Now suppose there is a map $g :\Omega \times \real^d \to \real^d$. For each $\omega \in \Omega$, one gets a different self-map $g(\omega, \cdot)$ on $\real^d$. If the choice of $\omega$ is random, then $g$ provides a parametric description of a stochastic process on $\real^d$. Now consider the following dynamical system on the product space $\Omega \times \real^d$:
\begin{equation*} 
	\begin{split}
		\omega_{n+1} &= \tau(\omega_n) \\
		y_{n+1} &= \mathfrak{f}\paran{ \omega_n, y_n }
	\end{split}
\end{equation*}
This is called a \emph{skew-product} system as the first variable evolves independently and continues to drive the dynamics in $\Omega$. Such skew-product systems  provide a universal description of discrete-time stochastic dynamics as a deterministic dynamical system \cite{Arnold1991random}. The stochasticity in the dynamics of the $y$ variable is interpreted to originate from the randomness of the initial state $\omega_0$. 

Next suppose that the deterministic map $f$ from Assumption \ref{A:f} corresponds to $\mathfrak{f}(\omega_0, \cdot)$ for some point $\omega_0 \in \Omega$. Now suppose that the function $\mathfrak{f}$ depends on a third parameter $\epsilon>0$ which represents a noise-bound. Thus $\mathfrak{f}$ may be denoted as $\mathfrak{f}_\epsilon$. Let $\alpha_\epsilon$ be an invariant measure for the system corresponding to $\mathfrak{f}_\epsilon$. Then $\proj_\Omega \alpha_\epsilon$ is an invariant measure for the $\Omega$-dynamics. Suppose that as $\epsilon\to 0^+$ the projections $\proj_\Omega \alpha_\epsilon$ converge weakly to the Dirac-delta measure $\delta_{\omega_0}$. Then $f$ is interpreted to be the \emph{zero-noise limit} of the parameterized family $\SetDef{\mathfrak{f}_\epsilon}{\epsilon>0}$ of stochastic processes. Following \cite{CowiesonYoung2005}, the ergodic system $(\Omega, \mu, f)$ is said to be \emph{stochastically stable} if for any parameterized family $\SetDef{\mathfrak{f}_\epsilon}{\epsilon>0}$, any choice of invariant measures $\alpha_eps
$, if $\proj_\Omega \alpha_\epsilon$ converges to $\delta_{\omega_0}$, then $\proj_{\real^d} \alpha_\epsilon$ must converge weakly to $\mu$. In other words if $f$ is the zero-noise limit of a parameterized family $\braces{\mathfrak{f}_\epsilon}$ of stochastic processes, then for any choices of stationary measures $\alpha_\epsilon$, their projection into $\real^d$ must necessarily converge to $\mu$.

The projection $\proj_\Omega \alpha_\epsilon$ characterizes the spread in the parameter space $\Omega$ and thus the spread in the uncertainty on the dynamics on $\Omega$. If an ergodic system is stochastically stable, and if it is represented by a stochastic process with a small spread in uncertainty around $f$, then any invariant measure of the stochastic process must be close to $\mu$. The concept of stochastic stability provides a rigorous platform on which to assess the visibility and stability of ergodic measures. Stochastic stability is hard to established in general, and has only been demonstrated 
called \emph{SRB}-systems \cite{SRB_young}. 

The following result from \cite{Das2024zero} summarizes the basic principle for approximating a deterministic dynamics by a Markov process :

\begin{lemma} [Step-skew approximation of dynamics] \label{lem:4}
	\cite{Das2024zero} Let Assumptions \ref{A:f} and \ref{A:cover} hold, and additionally suppose that the mesh size of $\calU $ is some $\delta>0$. 
	\begin{enumerate} [(i)]
		\item Let $\mu_\delta$ be any stationary measure of the step-skew system \eqref{eqn:Mrkv:1}. Then 
		\[ \support( \mu_\delta ) \subset B \paran{ X, \delta } , \quad X \subset B \paran{ \support( \mu_\delta ), 2\delta } . \]
		\item For any $N\in \num$ and $\eta \in (0,1)$, the partition may be chosen so that with a probability of at least $\eta$, an $N$-length trajectory of the stochastic dynamics \eqref{eqn:Mrkv:1} equals an $N$-length observed trajectory of $f$, namely $ \braces{ f^n(x_0) }_{n=0}^{N}$ for some $x_0$.
		\item The ergodic system $(X, f, \mu)$ is the zero noise limit of the stochastic dynamics \eqref{eqn:Mrkv:1}, as $\delta$ approaches zero.
	\end{enumerate}
\end{lemma} 

Lemma \ref{lem:4}~(ii) thus establishes any dynamical system universally as a zero noise limit of a stochastic dynamical system. The first claim implies that the support of the $\mu_\delta$ converges in Hausdorff metric to the targeted attractor $X$. This is one of our primary goals, an approximation of the invariant region which is being sampled by data. Corollary \ref{corr:1} is a summary of Lemma \ref{lem:4}~(ii).

The step-skew system \eqref{eqn:Mrkv:1} describes a random walk on $\real^d$. If all the cells of the partition are isomorphic then one can restrict the walk to a compact topological space. Most practical realizations of covers or partitions create simplexes or disks as the cells, thus satisfying the following simple assumption  :

\begin{Assumption} \label{A:disk}
	Each cell of the partition $\calU$ from Assumption \ref{A:cover} is topologically a $d$-disk $\calD^d$.
\end{Assumption} 

Assumption \ref{A:disk} implies that there are homeomorphisms :
\begin{equation} \label{eqn:A:disk}
	\begin{tikzcd} D^d \arrow[rr, "h_i", "\cong"'] && U_i \end{tikzcd} , \quad 1\leq i \leq m .
\end{equation}
These homeomorphisms \eqref{eqn:A:disk} allows the dynamics to be expressed as the step-skew system \eqref{eqn:Mrkv:1} in which the state space is $\calS \times D^d$. The transitions $\tilde{\phi}_{s_n\to s_{n+1}}$ can be described as
\[ \tilde{\phi}_{s_n\to s_{n+1}} : D^d \to D^d, \quad y_n \,\mapsto\, h_{s_n}^{-1} \circ \phi_{ s_n\to s_{n+1} } \circ h_{s_n} ( y_n ) .\]
The system \eqref{eqn:Mrkv:1} has two components -- an autonomous Markov process on the state space $\calS$, and a dynamics on the space $\calD^d$ which depends on the states $s_n$ and $s_{n+1}$.

Step-skew dynamical systems are an important class of dynamics. They have been used to demonstrate a variety of robust and non-intuitive behavior in dynamical systems \cite[e.g.]{GorodetskiEtAl1999, KleptsynNalskii2004, IlyashenkoNegut2010, DiazGelfertRams2011rich}. Lemma \ref{lem:4} presents their importance as a universal approximator for arbitrary dynamical systems. A step-skew system is a discrete time process, and its driving dynamics  is a Markov process. The rest of the paper presents how such a system can be approximated by a continuous time flow, in which the driving dynamics is deterministic, continuous-time, and drives the main system intermittently.

\section{Main results} \label{sec:main} 

\begin{figure}[!t]
	\includegraphics[width=0.95\linewidth, height=0.5\textheight, keepaspectratio]{\figs 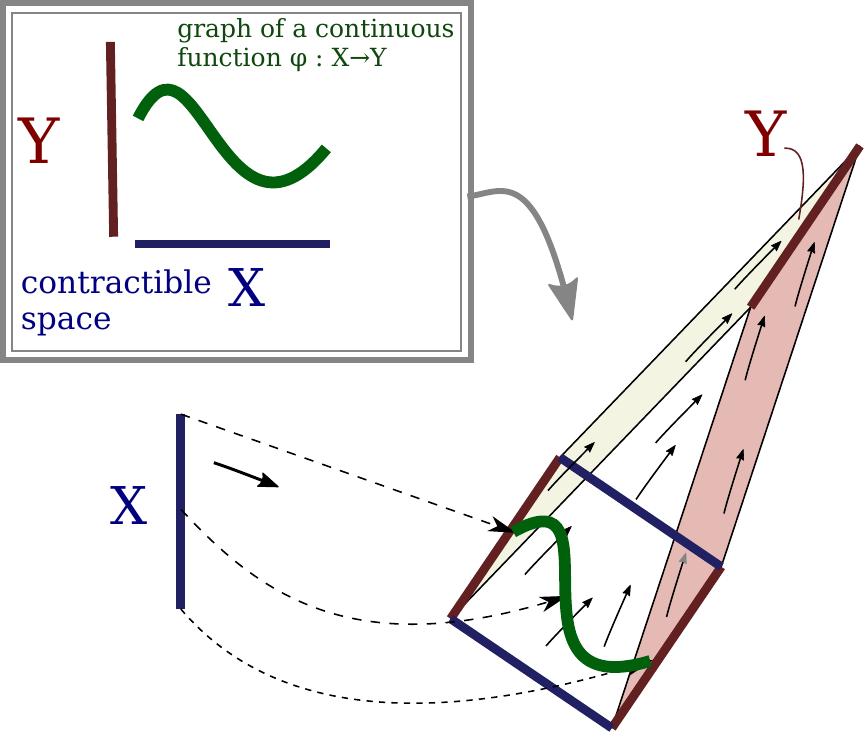}
	\caption{Converting a map into a flow along a pipe. The construction is based on a continuous map $f:X\to Y$ between topological spaces, with $X$ being contractible. The figure depicts a semi-flow from $X$ to $Y$ such that every point in $X$ is transported to the point $f(x)$ in $Y$. The semi-flow is based on the graph of $f$, which provides an embedding of $X$ in the product space $X\times Y$. The flow is divided into two stages. The first is a direct flow of each point $x$ into its corresponding point on the graph. The second stage is a semi-flow on the entire product space, which contracts it into $Y$. Let $H : [0,1] \times X \to X$ be the homotopy such that $H(0,\cdot) = \Id_{X}$ and $H(1,\cdot)$ is a constant map. This leads to the semi-flow $\Psi^t(x,y) = \paran{ H(t,x), y }$. }
	\label{fig:pipe}
\end{figure}

The Markov process \eqref{eqn:Mrkv:1} on the joint space induces a Markov process on $\tilde\calU := \cup_{i=1}^m U_i$. The Markov process on $\tilde\calU$ has the transition function : 
\begin{equation} \label{eqn:Mrkv:2}
	G : \tilde\calU \to \Prob \paran{ \tilde\calU }, \quad G(y) :=  \sum_{i : j\to i} \mathbb{P}_{i,j} \delta_{ \phi_{j\to i} (y) } , \quad \forall y \in V_j, \, 1\leq j \leq m .
\end{equation}
Equation \eqref{eqn:Mrkv:2} provides an equivalent description in terms of a Markov transition function. Since the measure $\mu$ is ergodic, the Markov transition on $\calS$ has a unique stationary measure $\nu$.

\paragraph{Sample paths} We wish to compare the performance of the step-skew product \eqref{eqn:Mrkv:1}, a stochastic process with that of the time-3 flow map perturbed pipe-flow, a deterministic process. These two dynamics which are fundamentally different may be compared by their space of \emph{sample-paths}. It is a standard practice to view a stochastic process as a random function on a time-domain. This was pioneered by Doob \cite{Doob1953book}, and subsequently a variety of conditions have been formulated \cite[e.g.]{GikhmanSkorokhod2004theory, Breiman1992prob} on the space of sample functions to characterize the resultant stochastic process. The Markov transition function \eqref{eqn:Mrkv:2} induces a discrete-time Markov process $\SetDef{X_n}{n\in\num_0}$ on $\real^d$. One can interpret this random process as a random variable 
\[ X : \Paths \times \num_0 \to \real^d  \]
with $\Paths$ being an abstract measurable space with a probability measure $\gamma$. This allows an interpretation of $X$ as the measurable map 
\[ X : \Paths \to \mathbb{F} \paran{ \num_0 ; \real^d } . \]
Thus a stochastic process may be interpreted as a random variable taking values in the space $\mathbb{F} \paran{ \num_0 ; \real^d }$ of functions from $\num_0$ to $\real^d$. In other words, a stochastic process is a random function from time $\num_0$ to state-space $\real^d$. This random variable $X$ projects the measure $\gamma$ into a measure $\law := X_* \gamma$ on $\mathbb{F} \paran{ \num_0 ; \real^d }$. This probability measure $\law$ will be called the \emph{law} of the stochastic process \cite{Pollard2012cnvrgnc}. 

\paragraph{Main result} We now state our main result, previously alluded to in Corollary \ref{corr:3} :

\begin{theorem} \label{thm:3}
	Suppose there is a step-skew product system as in \eqref{eqn:Mrkv:1} on the state space $\calS \times D^d$ where $\calS$ is a finite set and $\calD^d$. Let $\bar{\mu}$ be an invariant measure for the process. Then 
	\begin{enumerate} [(i)]
		\item there is a skew-product deterministic flow $\Psi^t_T$ on a topological space $\Omega \times \calX$ depending on a time-scaling parameter $T>0$,
		\item and a map $\phi : \calX \to D^d$,
	\end{enumerate}
	such that the sequences of measurements $\braces{ \phi \paran{ \Psi_T^{3n} x_0 } }_{n=0}^{\infty}$, emerging from various initial points $x_0\in \calX$, converge weakly in law to the series generated by \eqref{eqn:Mrkv:1} as $T$ goes to infinity.
\end{theorem} 

The flow promised by Theorem \ref{thm:3} will be called a \textit{perturbed pipe flow} (p.p.f.). There are three components to it. The first is the phase space $\calX$. This is a topological realization of the directed graph on the finite state-space $\calS$ which generates the Markov transition, as shown in Figure \ref{fig:network}. A directed graph comprises of nodes and edges. The edges will be realized as cylindrical \emph{pipes} and the nodes with multiple inputs and outputs will be realized as a special construction called a \emph{junction}. The $j$-th junction have one input window and multiple output faces. The number of output faces coincides with is the out-degree of vertex $j$ in the Markov transition graph. See Figure \ref{fig:cell} for an illustration. Overall the space $\calX$ remains embedded within an Euclidean space $\real^D$. The construction of junctions is described in Section \ref{sec:junction}. 

\begin{figure}[!t]
	\center
	\includegraphics[width=0.47\linewidth, height=0.5\textheight, keepaspectratio]{\figs 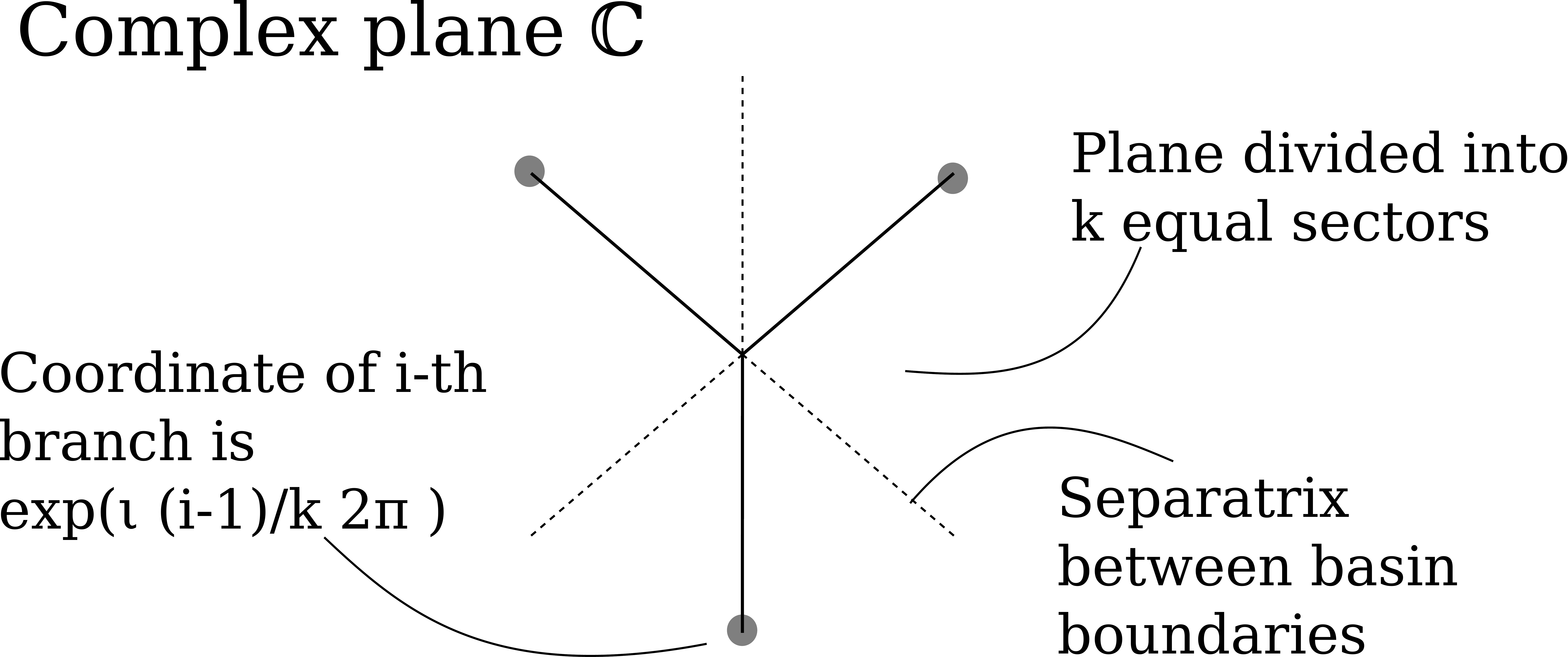}
	\includegraphics[width=0.47\linewidth, height=0.5\textheight, keepaspectratio]{\figs 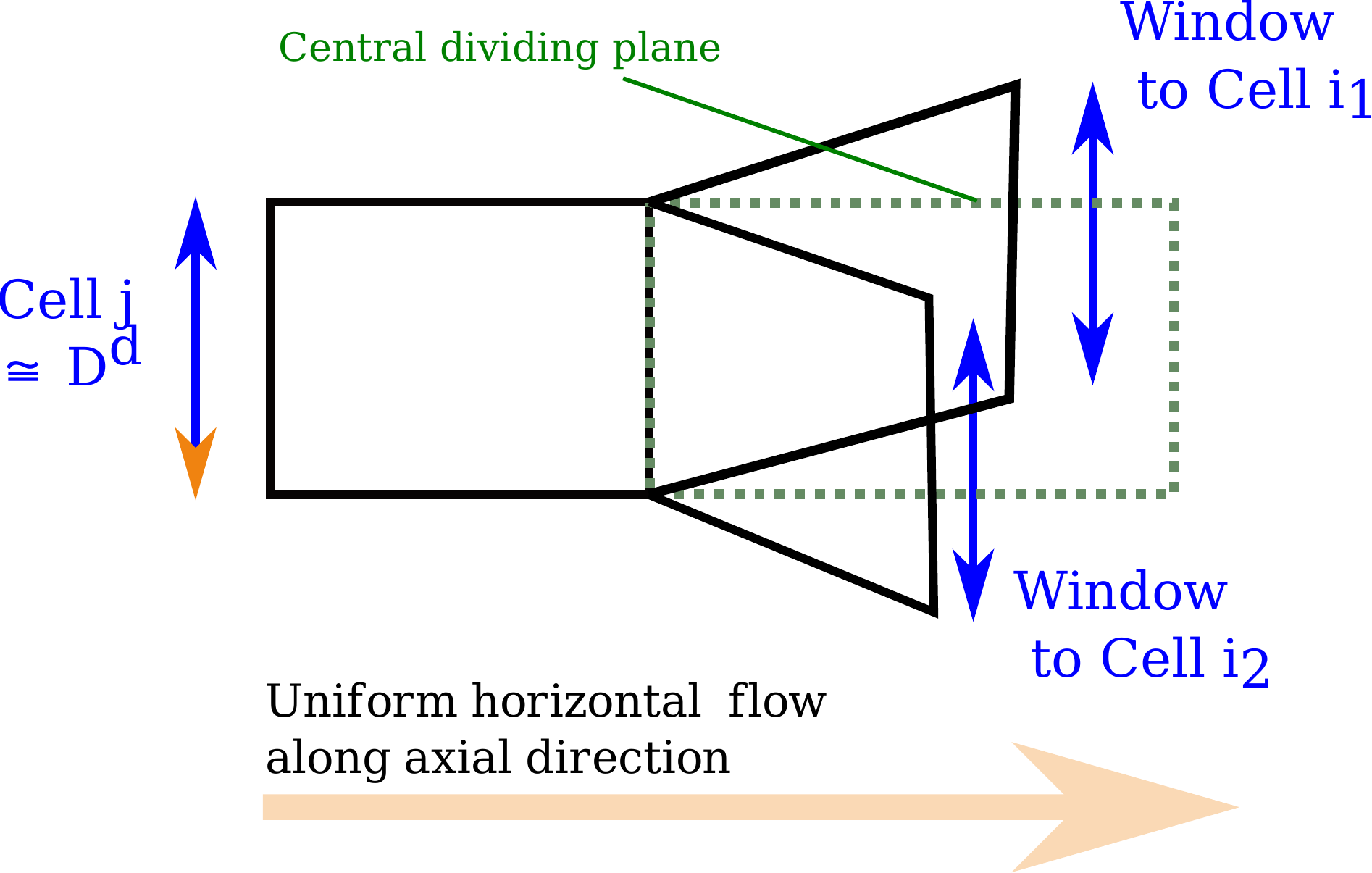}
	\caption{An frontal view along the axis of a junction. Figure \ref{fig:cell} presented a lateral view of a $k$ junction. To place the construction in $\real^{d+3}$ one starts by constructing the following branched 1-manifold in $\real^3$. The first branch representing the input channel, is a straight line segment from $(-1,0,0)$ to $(0,0,0)$. For each $1\leq i \leq k$, the $i$-th output channel is the straight line from $(0,0,0)$ to $(1, \cos \paran{ \frac{i-1}{k} 2\pi }, \sin \paran{ \frac{i-1}{k} 2\pi })$. Part (a) presents a frontal view of this branched 1-manifold. One next takes a Cartesian product with $\calD^d$ as shown in Figure \ref{fig:cell} to construct the junction. Every $k$-junction in $d$-dimensions constructed in Figure \ref{fig:cell} is provided a vector field. The vector field can be decomposed into two components - axial and lateral. The figure on the right describes the axial component of the flow. The diagram assumes $k=2$ for the sake of simplicity. The axial vector field is kept constant and equal to $1$ throughout the length of the junction. See Figure \ref{fig:cell_lateral} for a description of the lateral component of the vector field.}
	\label{fig:frontal}
\end{figure}

The second component of Theorem \ref{thm:3} is the flow $\Psi_T^t$. This flow will have distinctive characteristics in the pipes and junctions. Within each pipe corresponding to an edge $i\to j$, the flow is axial and laminar. The flow will be designed so that the net transport it causes equals the transition $\tilde{\phi}_{i\to j}$ as listed in \eqref{eqn:Mrkv:1}. Figure \ref{fig:pipe} presents an outline of this design. The more complicated part of the flow $\Psi_T^t$ is within each junction. Here the flow has both axial and lateral components, as illustrated in Figures \ref{fig:frontal} and \ref{fig:cell_lateral} respectively. The lateral component is driven by an external mixing dynamical system. Section \ref{sec:random} describes how this lateral component makes the flow behave like a pseudo-random number generator. As a result any point entering the $j$-th junction is directed towards one of its outgoing pipes with a probability distribution close to the vector $\transP(j)$. Thus overall the flow $\Psi_T^t$ is a deterministic, continuous-time dynamical system. We call this a {perturbed pipe flow} and provide the full details in \eqref{eqn:thm:3}.

The third component of Theorem \ref{thm:3} is the design of the external driving vector field, so that the sequence of switches made in successive junctions mimics the iterates of \eqref{eqn:Mrkv:1}. The flow $\Psi_T^t$ is designed so that at each of the instants $\SetDef{3n}{n=0,1,2,\ldots}$ the orbit is at one of the incoming faces of a junction. Moreover, each of these faces is homeomorphic to $\calD^d$. The purpose of the measurement $h$ is simply to relay the location of a point on a window of the junction into coordinates on $\calD^d$. Note that $h$ could take any form outside these sets without affecting the results.

The details of the construction are postponed to the later sections. We next establish the approximation properties of Theorem \ref{thm:3}.

\paragraph{Sequences} A p.p.f. is essentially a deterministic skew product system with two components --
\begin{enumerate} [(i)]
	\item A continuous time, deterministic, measure-preserving dynamical system $\Gamma_T^t : \Omega \to \Omega$ on some abstract probability space $(\Omega, \nu)$.
	\item A driven dynamical system on $\real^D$.
\end{enumerate}

The subset of interest to us in a compact, invariant subset of $\real^D$ which we denoted as $\calX$. The essential features of the flow are as follows :

\begin{lemma} \label{lem:pipe_shift}
	A perturbed pipe-flow (p.p.f) is a skew-product flow $\Psi_T^t$ driving a dynamical system on a topological space $\calX$ which can be partitioned into two components - pipes and junctions. Almost every orbit the flow satisfy the following --
	\begin{enumerate} [(i)]
		\item There is a finite collection $\calS$ of subsets transversal to the flow. These are called \textit{input windows}.
		\item For any initial condition of the driven variable on an input window, the orbit enters an input window every 3 time units.
		\item Let the itinerary created by the successive visits to windows be denoted as $\braces{s_n}_{n=0}^{\infty}$. Then the transitions $s_{n} \to s_{n+1}$ follow the paths of Markov transition matrix $\mathbb{P}$ on $\calS$.
		\item The time-3 samples bear the relation
		\begin{equation} \label{eqn:lem:pipe_shift}
			\Psi_T^{3(n+1)} ( \omega_0, x_0 ) = \phi_{s_n \to s_{n+1}} \paran{ \Psi_T^{3(n)} ( \omega_0, x_0 ) }, \quad \forall n\in \num_0 .
		\end{equation}
	\end{enumerate}
\end{lemma}

Lemma \ref{lem:pipe_shift} and in particular \eqref{eqn:lem:pipe_shift} establishes the connection between the skew-product dynamics \eqref{eqn:Mrkv:1} and the continuous time flow. The space $\calX$ is constructed so that the input windows serve as a physical representative of the states of $\calS$. The Markov transitions $\phi_{s_n \to s_{n+1}}$ are actuated by time-intervals of $3$ during which the orbit is transported by $\Psi_T$ from one window to another. At each junction, the orbit has several directions to choose from. Each direction corresponds to a transport to a different window. The choice of the direction depends on the driving dynamics $\Gamma^t_T$. Condition (iii) above states that the transitions conform to the transitions allowed by a prefixed Markov transition matrix $\mathbb{P}$. Condition (ii) guarantees that at teach time instant $3n$, we can keep track of the current cell / window $s$ instead of the state $\omega\in \Omega$. Thus a typical orbit of $\Psi_T^t$ may be represented by a sequence $\braces{ (s_n, y_n) }_{n=1}^{N}$. In fact the time-3 flow map of $\Psi_T^t$ takes the form of a discrete-time skew-product system
\begin{equation} \label{eqn:Mrkv:4}
	\begin{split}
		\omega_{n+1} &= \tau(\omega_n) \\
		y_{n+1} &= h \paran{ \omega_n, y_n } 
	\end{split}.
\end{equation}
We use to notation
\[ \braces{ (s_n, y_n) }_{n=1}^{N} \triangleright \; \eqref{eqn:Mrkv:4}.\]
to indicate that this sequence of pairs of states and $\real^D$-points are generated by such a skew-product system. In each of the skew product dynamics such as \eqref{eqn:Mrkv:1} or \eqref{eqn:Mrkv:4}, the driven dynamics is deterministic. Thus the distribution of the above sequences is dictated by the distribution of sequences of states from $\braces{1, \ldots, m}$. We study this next.

\paragraph{Conditional probabilities} Step-skew products are essentially a stochastic, discrete-state process driving a deterministic dynamics on some topological space $\calD$. Each trajectory of a step-skew system is a sequence of points in $\braces{1,\ldots, m} \times \calD$. When we consider the trajectories in the space of sample paths, their distribution is determined only by their component along $\braces{1,\ldots, m}$. This is made precise below : 

\begin{lemma} \label{lem:step_skew_shift}
	Let $\calS$ be a finite collection of states, $\mathbb{P}$ a Markov transition matrix on $\calS$, leading to a step-skew product system as in \eqref{eqn:Mrkv:1}. Let $\law$ be an stationary measure on the path space. Then
	\begin{enumerate} [(i)]
		\item The projection of $\law$ to the underlying symbolic space is a shift-invariant measure.
		\item The conditional probabilities follow the law
		\begin{equation} \label{eqn:step_skew_shift}
			\law \paran{ \braces{ (s_n, y_n) }_{n=1}^{N} \triangleright \; \eqref{eqn:Mrkv:1} \,|\, (s_0, y_0) } = \alpha \paran{ \Cylndr \paran{ s_0, \ldots, s_N } \,|\, s_0 } = \prod_{n=1}^{N} \mathbb{P}\paran{ s_{n-1}, s_{n} } ,
		\end{equation}
		where $\mathbb{P}$ is an $m\times m$ matrix of transition probabilities.
	\end{enumerate}
\end{lemma}

The proof of Lemma \ref{lem:step_skew_shift} follows directly from the fact that in a step-skew system, the iterates in the driven variable are deterministic functions of the the Markov state and the driven variable. A step-skew product can produce several sample paths from the same initial condition $(s_0, x_0)$. These sample paths are in one-to-one correspondence with the Markov walks starting from $s_0$. Thus all probability calculations originate in the finite-state Markov process. 

\paragraph{Proof of Theorem \ref{thm:3}} Lemma \ref{lem:step_skew_shift} applies to any step-skew dynamics, and thus to \eqref{eqn:Mrkv:1} and \eqref{eqn:Mrkv:4}. Let $\law$ be a shift-invariant probability measure of \eqref{eqn:Mrkv:1}, and $\tilde{\law}_T$ be a shift-invariant probability measure for the symbolic dynamics generated by \eqref{eqn:Mrkv:4}. In particular we get :
\[\begin{split}
	&\tilde{\law}_T \paran{ \braces{ \Psi_T^{3n} (s_0, y_0) }_{n=1}^{N} \,|\, (s_0, y_0) } = \tilde{\mu} \paran{ \Cylndr \paran{ s_0, \ldots, s_N } \triangleright \; \eqref{eqn:Mrkv:4} \,|\, s_0 } \quad \mbox{ by } \eqref{eqn:step_skew_shift}, \\
	&\quad \quad \quad \xrightarrow{T\to \infty} \prod_{n=1}^{N} \mathbb{P}_{ s_{n} , s_{n-1} }  \quad \mbox{ by Lemmas \ref{lem:uyg75}, \ref{lem:mixing:4} },\\
	&\quad \quad \quad = \mu \paran{ \Cylndr \paran{ s_0, \ldots, s_N } \,|\, s_0 } \\
	&\quad \quad \quad = \law \paran{ \braces{  (s_n, y_n) }_{n=1}^{N} \triangleright \; \eqref{eqn:Mrkv:1} \,|\, (s_0, y_0) } \quad \mbox{ by } \eqref{eqn:step_skew_shift} .
\end{split}\]
In conclusion we have proved that 
\begin{equation} \label{eqn:thm:3}
	\lim_{T\to\infty} \tilde{\law} \paran{ \braces{ \Psi_T^{3n} (s_0, y_0) }_{n=1}^{N} \triangleright \; \eqref{eqn:Mrkv:4} \,|\, (s_0, y_0) } = \law \paran{ \braces{  (s_n, y_n) }_{n=1}^{N} \triangleright \; \eqref{eqn:Mrkv:1} \,|\, (s_0, y_0) }.
\end{equation}
%
This type of convergence of  $\tilde{\law}$ to $\law$ is called \emph{weak conditional} convergence \cite{Sweeting1989cond}. Thus overall,  we establish that the deterministic process \eqref{eqn:Mrkv:4} has \emph{weak conditional convergence in law} to the step-skew product \eqref{eqn:Mrkv:1}. Equation \eqref{eqn:thm:3} formalizes the definition of weak conditional convergence in law. \qed

It still remains to establish the convergence claims upheld by Lemmas \ref{lem:uyg75}, \ref{lem:mixing:4}. The proof of these lemmas, as well as of Lemma \ref{lem:pipe_shift} will be presented as an outcome of the design details. This is done over the course of the next three sections.

\section{Random number generation} \label{sec:random} 

The purpose of this section is to show how chaotic systems behave as a pseudo-random generator. This behavior stems from a property called \emph{mixing}.

\paragraph{Mixing} Recall that a flow $(\Omega, \Gamma^t)$ is \emph{mixing} \cite[e.g.]{Halmos1944mixing, Nadkarni1998book, DasGiannakis_delay_2019} with respect to an invariant measure $\nu$ if for any two functions $\phi, \phi' \in L^2(\nu)$, one has the \emph{decay of correlations}  : 
\begin{equation} \label{eqn:mixing:1}
	\left\langle \phi, \phi'\circ \Gamma^t \right\rangle_{L^2(\nu)} = \int \phi \cdot \paran{ \phi' \circ \Gamma^t } d\nu \xrightarrow{t \to \infty} \norm{\phi}_{L^2(\nu)} \norm{\phi'}_{L^2(\nu)}. 
\end{equation}
In particular, if we choose $\phi, \phi'$ to be the indicator functions of two measurable sets $A, B$, then the rule of decay of correlations implies that 
\begin{equation} \label{eqn:mixing:2}
	\nu \SetDef{ x\in A }{ \Gamma^t x \in B } \xrightarrow{t \to \infty} \nu (A) \nu(B).
\end{equation}
Thus \eqref{eqn:mixing:2} indicates that in terms of correlations, the events $A, B$ appear independent after sufficient passage of time. This is due to innate nature of the flow to mix the points in the phase space. In other words :
\begin{equation} \label{eqn:mixing:3}
	\Condprob[\nu]{ x\in A }{ \Gamma^t x \in B } \xrightarrow{t \to \infty} \nu (A) .
\end{equation}
Equation \eqref{eqn:mixing:3} is the basis of the simple principle that a mixing but deterministic system can mimic a stochastic source. We however require a stronger form of mixing.

\paragraph{Multiple mixing} The ergodic system $\paran{\Omega, \Gamma^t, \nu}$ is said to be \emph{multiple-mixing of order $k$} \cite[e.g.]{ryzhikov1996stoch, Rue2023join} if for any $\phi_0, \ldots, \phi_k \in L^2(\nu)$ 
\begin{equation} \label{eqn:mixing:6}
	\lim_{ T_1, \ldots, T_{k-1}\to \infty } \nu \paran{ \phi_0 \cdot \prod_{i=1}^{k} \phi_i \circ \Gamma^{ -(T_1 + \ldots + T_i)} } = \prod_{i=0}^{k} \nu\paran{\phi_i} .
\end{equation}
A set theoretic formulation of \eqref{eqn:mixing:6} would mean that for any sets $A_0, \ldots, A_k$  :
\begin{equation} \label{eqn:mixing:7}
	\lim_{ T_1, \ldots, T_{k-1}\to \infty } \nu \SetDef{x\in \Omega}{ \Gamma^{T_1 + \ldots + T_i}(x) \in E_i ;\quad \forall 0\leq i\leq k } = \prod_{i=0}^{k} \nu\paran{E_i} .
\end{equation}
A sequence of times $t_1, \ldots, t_N$ will be called \emph{$T$-separated} if the difference between successive times on this sequence is $T$ or more. A consequence of \eqref{eqn:mixing:7} is that
\begin{equation} \label{eqn:mixing:8}
	\lim_{ \begin{array}{c} t_0, \ldots, t_{k} \;T \mbox{-separated} \\ T\to \infty \end{array} } \nu \SetDef{x\in \Omega}{ \Gamma^{t_i}(x) \in E_i ;\quad \forall 0\leq i\leq k } = \prod_{i=0}^{k} \nu\paran{E_i} .
\end{equation}
The limiting behavior of \eqref{eqn:mixing:8} will be the key principle used for random number generation. An ergodic system will be called \emph{multiple-mixing of all orders} or simply \emph{multiply mixing} if any of \eqref{eqn:mixing:6}, \eqref{eqn:mixing:7} or \eqref{eqn:mixing:8} holds for all $k$. Multiply mixing systems have been shown to exist in many natural flows arising in low dimensional manifolds \cite[e.g.]{Marcus1978horo, Savvidy2016spectrum, Savvidy1991monte}.

\paragraph{Pseudo-randomness} We have seen how an ordinary mixing system exhibits an independence of events by the phenomenon of decorrelation in \eqref{eqn:mixing:3}. A multiply mixing system allows this pseudo-independence to take place for any random sequence of events. This is made precise in the following lemma :  

\begin{lemma} \label{lem:mixing:4}
	Suppose that $\paran{ \Omega, \Gamma^t, \nu }$ is a deterministic, multiply  mixing system. Let $\calE$ be a finite collection of measurable subsets of $\Omega$ with nonzero $\nu$-measure. Then for every error bound $\epsilon>0$ there is a $T_0>0$ for which the following holds :  Take any $N\in \num$ and sequence of integers $n_1 < n_2 < \ldots < n_N$, and a sequence of events $E_1, \ldots, E_N $ drawn from $\mathcal{E}$. Then : 
	\begin{equation} \label{eqn:mixing:4}
		\forall T>T_0 \,:\quad \abs{ \prob[\nu]{ \Gamma^{T n_1} x_0 \in E_1, \, \ldots , \Gamma^{T n_N} x_0 \in E_N } - \prod_{n=1}^{N} \nu(E_n) } < \epsilon .
	\end{equation}
\end{lemma}

Lemma \ref{lem:mixing:4} is proved in Section \ref{sec:appendix}. Lemma \ref{lem:mixing:4} and \eqref{eqn:mixing:4} present a pseudo-randomness in the observations of a mixing system. If the observations are one among finitely many events, then there is a minimum wait time $T_0$ such that for any sequence of observations of the dynamics at intervals of length at least $T_0$, the events appear uncorrelated and independent. This pseudo-randomness is one of the main theoretical principles of our construction.  

We next describe a special choice of $(\Omega, \Gamma^t)$ that can help tailor the random number generation to any prior, discrete-event probability distribution within any desired error limit. 

\paragraph{Suspension flow} Let any discrete-time dynamical system $f:X\to X$ and a positive valued function $h:X\to (0, \infty)$ one can construct a continuous time system called a \emph{suspension flow} over the base map $f$ and with ceiling function $h$. The space for this flow is a quotient space
\[ X \times h / f := \SetDef{ (x,s) }{ x\in X, \, 0\leq s \leq h(x) } \,/\, \paran{ (x, h(x)) \sim (fx, 0) } .\]
The flow is given by 
\[\Gamma^t( (x,s) ) := 
\begin{cases}
	(x,s+t) & \mbox{ if } 0\leq s+t < h(x) , \\
	\Gamma^{t-h(x)}( (f(x), 0) ) & \mbox{ if } s+t \geq h(x)
\end{cases} , 
\quad \forall t\geq 0.\]
This continuous time-flow has the same degree of smoothness or continuity as the base map $f$. In fact the map $f$ turns out to be the return map of this flow to the section $\SetDef{(x,0)}{x\in X}$. If the ceiling function $h$ is constant valued, then the mixing properties of the $\Gamma^t$ is identical to that of $f$. For our purpose we chose $(f,X)$ to be the shift map on 2 symbols, and $h\equiv 1$. We denote the resultant shift map as $\generator$. This dynamics is useful as an explicit way to simulate random number generation, as we describe next.

\paragraph{Random number generation} Given the generator $\generator$ and a complex valued function $\zeta : \Omega_{\text{shift}} \to \cmplx$, we are interested in the ergodic sum : 
\begin{equation} \label{eqn:def:barzeta}
	\bar{\zeta}_T : \Omega_{\text{shift}} \to \cmplx, \quad \omega  \,\mapsto\, \int_0^T w\paran{ \frac{t}{T} } \exp \paran{ \iota \zeta \paran{ \generatorF^t \omega } } dt
\end{equation}
The function $w$ here is a weight function $w:[0,1]\to \real$ which is non-negative valued, integrable, and with integral $1$. The function $\bar{\zeta}_T ( \omega )$ is still a random variable, due to its dependence on $\omega$. As the parameter $T$ is increased the trajectory samples a greater portion of the space $\Omega$. Thus we should expect that as $T$ increases, $\bar{\zeta}_T$ reflects an averaged function, that has less fluctuations with respect to $\omega$. This is made precise below : 

\begin{lemma} \label{lem:d9dk4}
	Let $\paran{\Omega, \nu, \Gamma^t}$ be any ergodic system.
	Fix an $\epsilon>0$ a $k$-length vector $\beta$ with non-negative entries. Then there is an integrable function $\zeta : \Omega \to \cmplx$ and a $T_0>0$ such that for every $T>T_0$ the function $\bar{\zeta}_{T}$ from \eqref{eqn:def:barzeta} satisfies
	\begin{equation} \label{eqn:f3od0:1}
		\begin{split}
			&\paran{ \bar{\zeta}_{T*} \nu } \{(0,0)\} = 0. \\
			&\paran{ \bar{\zeta}_{T*} \nu } \braces{ z \mbox{ is a k-th root of unity } } = 0. \\
			&\paran{ \bar{\zeta}_{T*} \nu } \braces{ \frac{i-1}{k} 2\pi < \arg(z) < \frac{i}{k} 2\pi } \in \paran{ \beta_i-\epsilon , \beta_i + \epsilon }
		\end{split}
	\end{equation}
\end{lemma} 

Lemma \ref{lem:d9dk4} is a direct consequence of ergodic theory and is proved in Section \ref{sec:appendix}. We now describe the precise use of this lemma.

For each $1\leq i \leq k$ we use $\Theta_i$ to denote the angular sector $\SetDef{z\in \cmplx}{ \frac{i-1}{k} 2\pi < \arg(z) < \frac{i}{k} 2\pi }$. Lemma \ref{lem:d9dk4} thus says that the averaged output of the measurement $\zeta$ will fall with probability $1$ into one of these $k$-sectors, and the distribution will closely follow the vector $\beta$. Moreover, the flow $\generator$ is mixing since the base map is mixing. Therefore \eqref{eqn:mixing:4} applies, and combined with \eqref{eqn:f3od0:1} it implies that for any $\epsilon>0$ if $T$ is large enough, then : 
\begin{equation} \label{eqn:mixing:5}
	\forall \left\{ \begin{array}{c} n_1<n_2 < \ldots< n_X \in \num \\ i_1, \ldots, i_N \in \braces{1, \ldots, k} \end{array} \right.
	\,:\quad \prob{ \generatorF^{T n_1} x_0 \in \Theta_{i_1}, \, \ldots , \generatorF^{T n_N} x_0 \in \Theta_{i_N} } - \prod_{n=1}^{N} \beta_{i_n} < \epsilon .
\end{equation}
Equation \eqref{eqn:mixing:5} thus provides a concrete realization of a pseudo-random generator, in which the successive events are occur almost independently as one of $k$ sectors. This tool will be useful in the construction of our pipe-flow.

\section{Junctions} \label{sec:junction} 

\begin{figure}[!t]
	\center
	\includegraphics[width=0.95\linewidth, height=0.5\textheight, keepaspectratio]{\figs 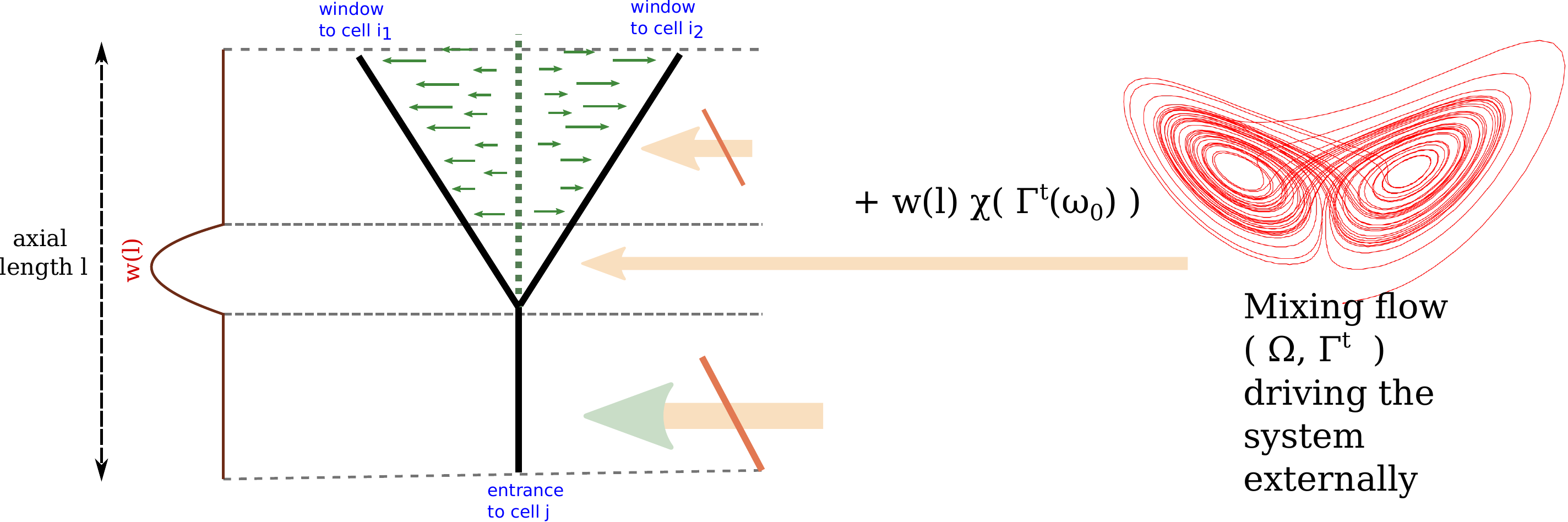}
	\caption{Lateral flow along a junction. The  $k$-junction in $d$-dimensions constructed in Figure \ref{fig:cell} is provided a vector field. The vector field can be decomposed into two components - axial and lateral. The figure presents a top-view of the junction, for the simple case when $k=2$. The axial coordinate is represented by a variable $l \in [0,2]$. the branching occurs at $l=1$. The lateral vector field is zero for $l\in [0,1.1]$. In the interval $[1.1,2]$ it is setup so that each of the branches are attracting sets, while the central axis remains neutral. Within the window $[1, 1.1]$ the junction receives a drift from an external mixing flow, shown in red. The input is weighted by the function $w$ from \eqref{eqn:def:weight}. Note that any trajectory under the combined action of the axial and lateral vector fields travels along the central axis till the branching point $l=1$. Within the window $1\leq l \leq 1.1$ it deviates to either of the branches. Due to the mixing nature of $(\Omega, \Gamma^t)$ and the uncertainty in its initial condition, this is a random event. Beyond the point $l=1.1$ the trajectory gets pulled to that branch in whose basin it lies. }
	\label{fig:cell_lateral}
\end{figure}

Junctions are continuous-time deterministic realizations of the various transitions occurring for the states $s$ in \eqref{eqn:Mrkv:1}, along with their transition probability measures $\transP(s)$. Each transition in a finite state Markov chain is a discrete-time event on discrete space. The probability vector $\transP(s)$ governs a random switch from current state $s$ into any one of the $m$ states of the Markov process. However, the number of output states with non-zero transition probability may be a number $k$ less than $m$. The junction that realizes such a transition to $k$ states is called a \emph{$k$-junction}. Figure \ref{fig:cell} presents the construction of a typical cell. Figure \ref{fig:frontal} presents how the $d$-dimensional cell can be embedded in $\real^{d+3}$. Their construction can be made precise by two introducing two notions - \emph{partial} semi-flows and \emph{axial, partial} semi-flows.

\paragraph{Partial semi-flows} A continuous \emph{semi-flow} on a topological space $\calX$ is an action of the semigroup $\paran{ \real^+_0, + }$ on $\calX$. It is semi-group homomorphism $\Phi$ from the additive semi-group $\paran{ \real^+_0, + }$ into the semi-group of endomorphisms of $\calX$. Thus for each $t\geq 0$ there is a continuous map $\Phi^t : \calX \to \calX$ such that for every $s+t\geq 0$, the rule $\Phi^{s+t} = \Phi^s \circ \Phi^t$ holds. A relaxation of these strict rules of a semi-flow is a \emph{partial semi-flow}. A partial semi-flow on $\calX$ consists of 
\begin{enumerate}
	\item a continuous function $\calT : \calX \to [0,\infty)$, to be interpreted as an \emph{exit-time}. The function $\calT$ assigns to every point $x$ the maximum time up to which there exists a path from $x$. .
	\item A subset $\tilde{\calX}$ of the product space $\calX \times [0,\infty)$ such that for each $x\in \calX$, the $x$-section of $\tilde{\calX}$ is $\{x\} \times [0, \calT(x)]$.
	\item A map $\Phi : \tilde{\calX} \to \calX$ such that 
	\begin{equation} \label{eqn:def:partial_flow:1}
		\forall (x,s) \in \tilde{\calX}, \quad \calT \paran{ \Phi(x,s) } = 0.
	\end{equation}
	The point $(x,s)$ indicates a point which is time $s$ ahead along the flow from $x\in \calX$. This address of $(x,s)$ on $\calX$ is given by $\Phi(x,s)$.
	The level set $\braces{ \calT = 0 }$ may be interpreted as the exit points of the set $\calX$. Thus $\Phi$ maps every point $a$ along with its exit time $\calT(a)$ into an exit point.
	\item Finally the following algebraic condition is met : 
	\begin{equation} \label{eqn:def:partial_flow:2}
		\forall (x,s) \in \tilde{\calX} , \quad  \calT \circ \Phi(x,s) = \calT(x) - s.
	\end{equation}
	The means that exit time of $\Phi(x,s)$, which is time $s$ ahead of $x$, is $s$ less than the exit time of $x$.
\end{enumerate}

\begin{figure} [!t]
	\includegraphics[width=0.48\linewidth, height=0.5\textheight, keepaspectratio]{\figs 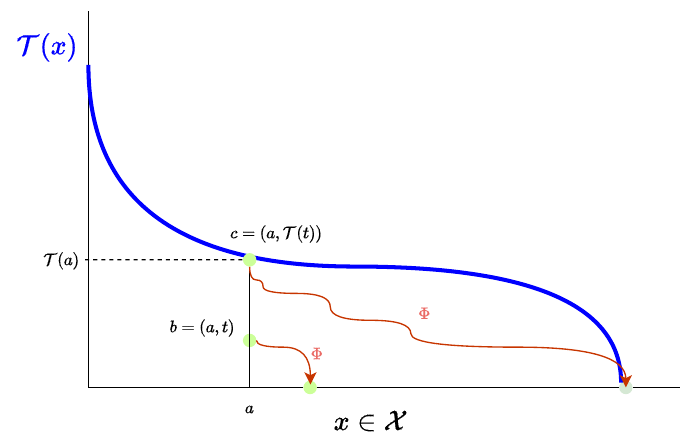}
	\includegraphics[width=0.48\linewidth, height=0.5\textheight, keepaspectratio]{\figs 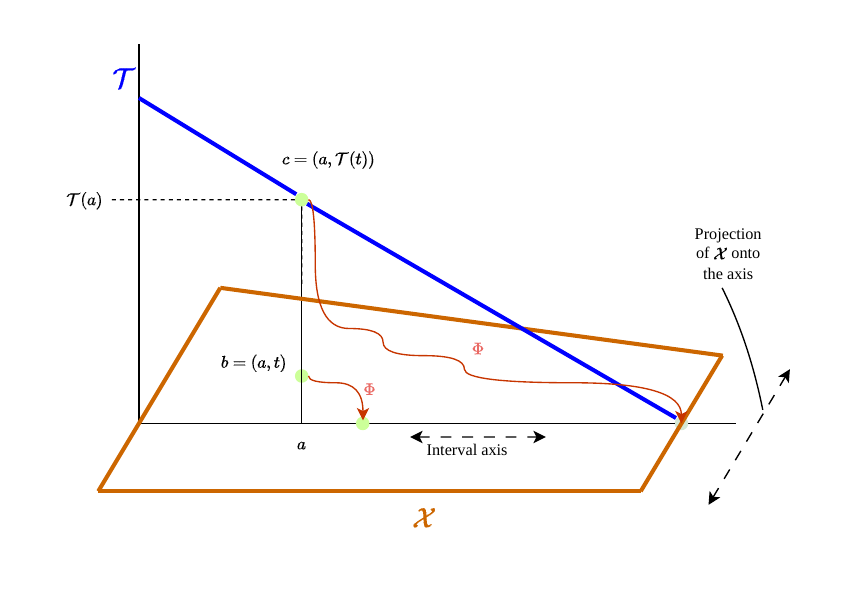}
	\caption{A partial semi-flow (left) and an axial semi-flow (right). These serve as templates for the type of constructions described in the paper. A partial semi-flow generalizes the notion of a semi-flow by allowing trajectories to be finite-time instead of extending indefinitely into the future. Each point $x$ on the phase space $\calX$ has an exit time $\calT$. The points on $\calX$ for which $\calT \equiv 0$ are called the exit points. The flow $\Phi$ maps every point $a$ along with its exit time $\calT(a)$ into an exit point. Given any time $s$ less than $\calT(a)$, $\Phi$ maps $(a,s)$ into a point at time $s$ further along the orbit of $a$. The space $\calX$ can be any topological space, and has been drawn as a line in the left figure for simplicity. When $\calX$ is indeed an interval and the exit time varies linearly, one gets an interval semi-flow \eqref{eqn:def:intrvl_part_flow}. An axial semi-flow is a semi-flow which can be projected into the interval semi-flow. Such a flow has been depicted on the right hand diagram. }
	\label{fig:partial_semiflow}
\end{figure}

See Figure \ref{fig:partial_semiflow} for an illustration. An usual semi-flow is a partial semi-flow, where $\Phi$ is simply the flow-map, and $\calT \equiv \infty$. Thus every initial state has an infinitely long orbit. In a general partial semi-flow, $\calT$ can be finite valued, meaning that infinitely long orbits do not exist. A simple example of a partial semi-flow which is not a semi-flow, is on the interval $[0,L]$. The exit-time and flow-map are respectively 
\begin{equation} \label{eqn:def:intrvl_part_flow}
	\begin{split}
		& \calT_{uni,L} : [0,L] \to [0 , L], \quad \calT_{uni,L}(l) := L-l. \\
		& \Phi_{uni,L} : (l, s) \mapsto (l+s, 0)
	\end{split}
\end{equation}
This concept of a partial semi-flow makes it easier to interpret the vector fields that we impart to the topological objects we construct. These objects -- pipes and junctions, come with a natural notion of a central axis. The flow that we impart to these objects have a uniform motion along their axial directions. This is made precise in our next definition.

\paragraph{Axial flows} An \emph{axial partial semi-flow} or simply an \emph{axial flow} is a partial semi-flow $\paran{\calX, \calT, \Phi}$ along with
\begin{enumerate}
	\item a continuous surjective map $\pi : \calX \to [0,L]$;
	\item for every $x\in \calX$, $\calT(x) = \calT_{uni,L} \circ \pi(x)$, which equals $L - \pi(x)$ by \eqref{eqn:def:intrvl_part_flow}; 
	\item and finally the identity :
	\[ \forall (x,s) \in \tilde{\calX}, \quad \pi\circ \Phi(x,s) = \pi(x)+s . \]
\end{enumerate}
Thus axial flows are partial semi-flows, in which there is the notion of an axis borne by a projection function $\pi$, and the projection creates a commutation between the uniform semi-flow on an interval and the original semi-flow. Due to this analogy, we shall call the subsets $\pi^{-1}(0)$ and $\pi^{-1}(L)$ the entry and exit faces of the axial flow. We call the quantity $L$ the \emph{length} of this axial flow. We next review some basic concepts from ergodic theory, the measure theoretic aspects of dynamical systems.

The flow through the junction is closely tied to its design. The $k$-junction has an axial direction, and a complimentary set of directions, collectively to be called the \emph{lateral} direction. This axis becomes easier to visualize when $k=2$, and is shown in Figure \ref{fig:frontal}~(b). The vector field imposed on the $k$-junction has axial and lateral components. The axial component is uniform and equal to 1, as shown in the figure. Thus any point starting from the initial window will travel at a uniform speed down the length of the junction. When its trajectory comes to the section where all the branches are glued, it receives a vector field perturbation in the lateral direction as well. As a result it makes a switch to one of the $k$ branches. The lateral vector field is zero all across the junction, except in the window $[1,1.1]$ as shown in Figure \ref{fig:cell_lateral}. The weight-function used in the latter figure is given by the formula :  
\begin{equation} \label{eqn:def:weight}
	w : \real \to \real, \, w(l) := 
	\begin{cases}
		(l-1)*(1.1-l) & \mbox{ if } 1 \leq l \leq 1.1 \\
		0 & \mbox{ otherwise }
	\end{cases}
\end{equation}
The choice of the branch depends upon the input it receives from the external flow $(\Omega, \Gamma^t)$ during the time interval $[1,1.1]$. The purpose of $w$ is to make the net vector field continuous over the junction. By making this external source $(\Omega, \Gamma^t)$ a mixing system, we can make this branching a random event. This is where we use the pseudo-randomness results of Section \ref{sec:random}. The conditions in \eqref{eqn:f3od0:1} are directly applicable to the layout of the exit windows of the junction, as displayed in Figure \ref{fig:frontal}~(a). 

\paragraph{Switchings} We set our excitation function $\chi : \Omega_{ \text{susp} } \to \real^{d+3}$ to be the function whose first coordinate and last $d$ coordinates are zero, and the second and third coordinates are $\exp\paran{ \iota \zeta (\omega) }$. The net lateral displacement within the axial window $[1, 1.1]$ can be calculated as follows : 
\[\begin{split}
	&\mbox{Total lateral displacement } \\
	&\quad \quad =  \proj_{2,3} \int_{1}^{1.1} w(l) \chi \paran{ \Gamma^{Tl} \omega_0 } dl  = \int_{1}^{1.1} w(l) \proj_{2,3} \chi \paran{ \Gamma^{Tl} \omega_0 } dl , \\
	&\quad \quad = \int_{1}^{1.1} w(l) \exp \paran{ \iota \zeta \paran{ \Gamma^{Tl} \omega_0 } } dl , \; \mbox{ take } l = 1 + s/T , \\
	&\quad \quad =  \frac{1}{T} \int_{0}^{0.1T} w \paran{ \frac{s}{T} } \exp \paran{ \iota \zeta \paran{ \Gamma^{s} \paran{ \Gamma^1(\omega_0) } } } ds  ,\\
	&\quad \quad = 0.1 \bar{\zeta}_{0.1T} \paran{ \Gamma^T(\omega_0) } \, \quad \mbox{ by \eqref{eqn:def:barzeta}}.
\end{split}\]
The function $\bar{\zeta}_{0.1T}$ in \eqref{eqn:def:barzeta} thus describes the total effect of the axial component of the vector field that the junction receives as input from the external ergodic flow. The value of $\bar{\zeta}_{0.1T}(\omega)$ decides which sector of the complex plane the an initial point starting at the origin is driven towards. Each such sector is the basin of attraction of one of the branches. Thus the value of $\bar{\zeta}_{0.1T}(\omega)$ decides which exit branch the trajectory is pull towards. The uncertainty in the location of $\omega$ within $\Omega$ is exactly the same as that of $\Gamma^T(\omega_0)$. This makes the switching a probabilistic event. 

Next we argue that in spite of a junction being driven by a completely deterministic external source, two successive selection of gates are almost-independent. The switch made in a junction depends implicitly on the initial state $\omega_0$ of the driver $\paran{ \Omega, \Gamma^t }$. Note that the junction receives an input from the driver only during the passage of the flow between the axial coordinates $[1,1.1]$, as indicated in Figure \ref{fig:frontal}~(b). The direction towards which the flow is nudged depends on the sum of the input it receives from the driver, via the function $\chi$, and during this period. Recall that the $\bar{\zeta}_{0.01T}$ can be designed so as to create the configuration in \eqref{eqn:f3od0:1}. Thus Lemma \ref{lem:d9dk4} guarantees the following : 

\begin{lemma} \label{lem:uyg75} 
	Consider a $k$-junction with transition probabilities $\beta = \paran{ \beta_1, \ldots, \beta_k }$,   the following holds :
	\begin{enumerate} [(i)]
		\item Consider any point $p$ at the entry window of a $d$-dimensional $k$-junction. Let its coordinate correspond to a point $x\in D^d$. Then the coordinate of the trajectory of $p$ continues to remain $x$ along $\calD^d$.
	\end{enumerate}
	Now consider any continuous-time ergodic flow $\paran{ \Omega, \Gamma^t, \nu }$, and an error tolerance $\epsilon>0$. Then there is a vector field $V$ as in \eqref{eqn:pipe_flow} such that for $T$ large enough,
	\begin{enumerate}[(i), resume]
		\item with probability $1$ the trajectory of $p$ exits through one of the $k$ exit-windows.
		\item The probability of exiting through the $i$-th window is within $\epsilon$ error of $\beta_i$.
	\end{enumerate}
\end{lemma} 

Thus overall we have a setup in which there is a topological object with a copy of $\calD^d$ at one end, and $k$-copies of $\calD^d$ at the other end. The former is called an entry window, and the latter are called the exit windows. Any trajectory travels uni-directionally and uniformly along the axis of this object and exits through one of the $k$ exit-windows. This event is actuated by an intermittently acting vector field which depends on the initial state in an ergodic system $(\Omega, \Gamma^t)$. The uncertainty in the initial state makes the choice of the exit window a probabilistic event. We have described how to design the intermittent vector field so that the probabilities of these events are within an $\epsilon$-tolerance of a prescribed vector of probabilities. This idea of linking a stochastic system with a skew-product system has been used in the reverse direction before \cite{MelbourneNicol2005almost, BurovEtAl2021kernel, MelbourneStuart2011skew}, by establishing a stochastic differential equation to be the limit of skew-product systems .

To approximate the step skew system \eqref{eqn:Mrkv:1} we shall need one junction each for each of the $m$ states. If we assume that the driving dynamics is mixing, then we can also establish that the successive switchings in different junctions are quasi-independent. We next describe how the junctions are connected together by pipes, and how the partial semi-flows on each of these components join together to form one single flow.

\section{Pipes} \label{sec:network} 

\paragraph{Pipes} The second component of our construction is a pipe. A pipe is a realization of any map $\phi :X\to Y$ as an axial flow in a topological space, whose initial states are in $X$ and final states in $Y$. Figure \ref{fig:pipe} presents the construction of such a pipe. Topologically, a pipe is created by gluing $X \times I$ with $X\times Y \times I$. The gluing occurs along the graph of the function $\phi$, as shown in the figure. This space is assigned a semi-flow as indicated by the arrows. The axial component of the semi-flow has a constant speed of $2$, so that the initial window gets transported to $Y$ in time $1$. 

In our case, for each transition $j\to i$, we join the $j$ and $i$-th cell with the pipe in which $X = Y = D^d$ and $g = \phi_{j\to i}$. This completes the construction of all components of the pipe-flow. We have provided separate constructions of the junctions and pipes. Each component is an axial flow. The axial flows transport points from their entry window into their exit face. We have thus built the following abstract setup : 

\paragraph{Network of axial flows} Let $\calG$ be a finite directed graph on $m$ vertices. Further suppose that for each $1 \leq j \leq m$ :
\begin{enumerate} [(i)]
	\item the vertex $j$ corresponds to a topological space $\calJ_j$ and an axial flow $\paran{ \calJ_j, \calT_j, \Psi_i^t, \pi_j, L_j }$.
	\item the exit window of $\paran{ \calJ_j, \calT_j, \Psi_i^t, \pi_j, L_j }$ has as many connected components as the out-degree of $j$;
	\item for every edge $j\to i$ from $j$ in $\calG$, we call the corresponding component of the exit window to be the window with index $i$;
	\item each edge $j\to i$ corresponds to a topological space $\Pipe_{j\to i}$  and an axial flow $\paran{ \calJ_{j\to i}, \calT_{j\to i}, \Psi^t_{j\to i}, \pi_j, L_{j\to i} }$; 
	\item the entry window of $\Pipe_{j\to i}$ is homeomorphic to the index-$i$ exit window of $J_j$, and whose exit window is isomorphic to the entry window of $J_i$;
	\item the axial flow on the $j$-th vertex is denoted as $\Psi_{j}^t$, and the axial flow on the edge $j\to i$ is denoted as $\Psi_{j\to i}^t$.
\end{enumerate}
The following lemma describes how these can be joined together to give a single consistent flow.

\begin{lemma} \label{lem:join_flows}
	Suppose one has the network of axial flows as described above. Then one can create a topological space by glueing the index-$i$ exit window of junction $j$ to the cylinder $j\to i$, and the exit window of the cylinder $j\to i$ to the entry-window of junction $i$. Let $a$ represent any component of the network, either a vertex or an edge. This space has a semi-flow defined recursively as 
	\begin{equation} \label{eqn:def:pipe_flow}
		\Phi^t (x) := 
		\begin{cases}
			\Phi^{t - 1 + \pi_a(x)} \paran{ \Psi_i^{ 1-\pi_a(x) } (x) } & \mbox{ if } t \geq L_a - \pi_a(x) ,\\
			\Psi_i^{ t} (x) & \mbox{ if } t \leq L_a - \pi_a(x) ,\\
		\end{cases}
	\end{equation}
	for every $x$ in component-$a$. An analogous result holds if the axial flows along the components are non-autonomous.
\end{lemma}

The flow \eqref{eqn:def:pipe_flow} is our desired realization of the step-skew system \eqref{eqn:Mrkv:1}. Figure \ref{fig:network} depicts how a step-skew system such as \eqref{eqn:Mrkv:1} is converted into a pipe flow - which is a $d+1$ dimensional branched manifold. 
\begin{enumerate}
	\item Each of the states $\braces{1, \ldots, m}$ of the Markov transition become junctions. If state $s\in \braces{1, \ldots, m}$ has $k$ outgoing states, then it becomes a $k$-junction.
	\item Each transition $j\to i$ is realized via a cylinder $\Pipe_{j\to i}$. It is provided the flow $\Psi^t_{j\to i}$ from \eqref{eqn:pipe_flow} .
	\item The entry window of the cylinder $\Pipe_{j\to i}$ is identified with the $i$-th exit window of junction $j$.
	\item The exit window of the cylinder $\Pipe_{j\to i}$ is identified with the entry window of junction $i$.
\end{enumerate}

The perturbed pipe-flow so constructed proves the statement of Corollary \ref{corr:2}.	We next discuss the approximation properties of the flow $\Psi_T^t$. 

\section{Conclusions} \label{sec:conclus} 

We have seen the construction of the perturbed pipe-flow in Sections \ref{sec:junction} and \ref{sec:network}. The partial flows through the junctions, and pipes are joined together with the help of Lemma \ref{lem:join_flows} to create a continuous flow for the entire topological space. The flow has constant speed of $1$ through the axial directions of each junction and pipe. Any trajectory enters a gate at fixed time intervals of $3$. The trajectory takes time $2$ to traverse through a junction. During this time period it makes a switch to one of the exit channels of the junction. The switching is actuated by the effect of the external driving $\Psi^t_{\psi}$. The special nature of step-skew products allow a simplified calculation of conditional probabilities in the path-space, as laid out in \eqref{eqn:step_skew_shift}. This leads to the weak conditional convergence in law, claimed in Theorem \ref{thm:3}.

We have thus presented three types of dynamical systems in Figure \ref{fig:outline}, and described how they perform approximations of each other. Any deterministic, discrete-time dynamics may be approximated by a step-skew dynamical system, in which a finite state Markov process drives a deterministic dynamical system. The approximation is in terms of both the invariant set as well as path space. Given any dynamics of this second kind, one can perform yet another approximation by a continuous time, deterministic skew-product system. The approximation is also in the conditional probabilities on the space of sample paths.

As mentioned before, the idea that a deterministic  chaotic system might mimic a stochastic process has been verified theoretically from multiple angles.  While non-constructive proofs have provided decisive results on the limits of learning and reconstruction, an explicit construction provides produces a tangible object of study. The article aims to provide a blueprint for studying chaotic dynamical systems in the same vein as Geometric Lorenz attractor \cite[e.g.]{GuckenheimerWilliams1979strctrl, ShiWang2025measures}, Logistic map \cite[e.g.]{Jakobson1981acim}, Blenders \cite[e.g.]{BonattiDiaz1996blender, BonattiDiaz2008blender}, Multi-chaos \cite[e.g.]{DasJim2017chaos} and Chua-circuits \cite[e.g.]{Matsumoto1984chua, Chua1986chua, Chua1995chua}. 

Two benefits may be reaped from an explicit construction. Firstly it provides a tool for computation and simulation. Secondly it provides a new mechanism for the emergence of mixing while preserving a specific topology. This mechanism may be tweaked in multiple ways to attain further properties associated a dynamical system. For example one may attempt to control the Lyapunov exponents by controlling the rate of contraction or expansion in a neighborhood of the invariant set. One can also try to achieve specific mixing and multiple-mixing rates. All of these provide interesting topics for deeper study.

We end with a particular implementation challenge that would make the reconstruction closer to the original dynamical system. Note that the dimension of the pipe-flow does not depend on the choice of partition but rather, on the dimension $d$ of the  phase space as specified in Assumption \ref{A:f}. A finer partition would lead to a closer approximation of the dynamics of $f$ as discussed before. If we take a progressively finer sequence of measurable partitions of a neighborhood of the attractor $X$ by repeated sub-division, then we get a sequence of Markov walks $\braces{ \paran{ \calS_n, \mathbb{P}_n } }_{n=1}^{\infty}$ in which $\paran{ \calS_n, \mathbb{P}_n } $ is a coarse-grained version of $\paran{ \calS_{n+1}, \mathbb{P}_{n+1} } $. The resulting pipe-flow realizations $\Gamma_{T,n}^t$ will still have the same dimension, but with more junctions and branchings. In fact the growth in the number of junctions is exponential with an asymptotic rate equal to the entropy of the system \cite[e.g.]{Halmos1956, Das2023CatEntropy}. One can choose an ambient Euclidean space $\real^D$ of sufficiently high dimension $D$ so that each of the the $\Gamma_{T,n}^t$-s can be embedded within $\real^D$. One can now ask whether the embedding can be chosen so that adjacency information from the original space $\real^d$ is preserved. This would lead to a more complete recreation of the dynamical system $(f,X)$ in $\real^D$. Moreover, one could study spatial perturbations of the reconstructed system to get meaningful insights into the stability of the original dynamics.
\section{Appendix} \label{sec:appendix}

\paragraph{Proof of Lemma \ref{lem:mixing:4}} Fix an $N\in \num$. Define events
\[E_i := \SetDef{x\in \Omega}{ \Gamma^{Ti}x \in A_i } , \quad 0\leq i \leq N .\]
Then note that there $T(N, \epsilon)>0$ such that for every $T>T(N, \epsilon)$, we have
\[ \abs{ \Condprob[\nu]{ E_{0} \cap E_1 \cap \cdots \cap E_m }{ E_{0} \cap E_1 \cap \cdots \cap E_{m-1} } - \nu(E_m) } < \epsilon , \quad 1\leq  . \]
Note that the probability $\nu(E_m)$ is just $\nu(A_m)$. Since the numbers $\nu(A_0), \ldots, \nu(A_N)$ are finite and non-zero, we can refine the choice of $T(N, \epsilon)>0$ such that
\[ \nu(A_m)^{-1} \Condprob[\nu]{ E_{0} \cap E_1 \cap \cdots \cap E_m }{ E_{0} \cap E_1 \cap \cdots \cap E_{m-1} } \in (1-\epsilon, 1+\epsilon) , \quad 1\leq m \leq N . \]
Then we have
\[\begin{split}
	& \Condprob[\nu]{ \Gamma^{Ti}x \in A_i , \; 0\leq i\leq N }{ x\in A_0 } = \Condprob[\nu]{ E_{0} \cap E_1 \cap \cdots \cap E_N }{ E_{0} } = \prod_{m=1}^{N} \Condprob[\nu]{ E_{0} \cap E_1 \cap \cdots \cap E_m }{ E_{0} \cap E_1 \cap \cdots \cap E_{m-1} }.
\end{split}\]
Thus
\[\Condprob[\nu]{ \Gamma^{Ti}x \in A_i , \; 0\leq i\leq N }{ x\in A_0 } = \prod_{i=1}^{N} \nu(A_i) \prod_{i=1}^{N} \gamma_i, \quad \gamma_i \in (1-\epsilon, 1+\epsilon) . \]
Thus w can set bounds
\[ (1-\epsilon)^N \leq \SqBrack{\prod_{i=1}^{N} \nu(A_i)}^{-1} \Condprob[\nu]{ \Gamma^{Ti}x \in A_i , \; 0\leq i\leq N }{ x\in A_0 } \leq (1+\epsilon)^N . \]
This completes the proof of Lemma \ref{lem:mixing:4}. \qed

\paragraph{Proof of Lemma \ref{lem:d9dk4}} We shall prove a stronger version of the Lemma. We shall show that given any $L^2(\nu)$ function $\psi : \Omega \to \cmplx$, one can find a $\zeta$ such that $\bar{\zeta}$ is arbitrarily close to $\psi$. If this is true, then the projection of $\nu$ under $\bar{\zeta}$ will also be arbitrarily close to the projection under $\psi$. The function $\psi$ can be chosen to achieve a probability measure on $\cmplx$ which is concentrated around the $k$ roots of unity, with weights $\beta_1, \ldots, \beta_k$ as described. This would achieve the approximation result.

The flow $\Gamma^t$ induces the following transformation on the function space $L^2(\nu)$  :
\[ \phi \mapsto \phi \circ \Gamma^t . \]
This transformation is a unitary group of operators, called the \emph{Koopman group}. For each time $t$, the transformation is denoted by $U^t$. Associated to such a group of unitary operators is a projection valued measure $E$. This is a function which maps each Borel subset on the unit circle $S^1$ in the complex plane, into a projection operator on $L^2(\nu)$. This assignments satisfies the usual axioms of a measure, such as additivity and zero value on the null set. Aided with this spectral measure one can write the unitary operator as an integral of an operator valued measure
\[U^t = \int_{S^1} e^{\iota \theta t} dE(\theta) . \]
This notation simplifies the notation of operator algebra. In particular we can write
\[\int_0^T w\paran{ \frac{t}{T} }U^t dt = \int_{S^1} \SqBrack{ \int_0^{T} w\paran{ \frac{t}{T} } e^{\iota \theta t} dt } dE(\theta) . \]
Now since for each $t$ $U^t$ is a unitary operator, $U^t$ is invertible. As a result, the operator expressed as the integral above has dense range for every $T>0$. Thus given any error bound $\epsilon>0$ and a function $\psi\in L^2(\nu)$ there is a function $\zeta\in L^2(\nu)$ such that 
\[ \norm{ \int_0^T w\paran{ \frac{t}{T} }U^t \zeta dt - \psi }_{L^2(\nu)} < \epsilon . \]
Now note that the function $\int_0^T w\paran{ \frac{t}{T} }U^t \zeta dt$ is precisely the function $\bar{\zeta}$ in the Lemma. This completes the proof of Lemma \ref{lem:d9dk4}. \qed


\end{document}